\documentclass[final, leqno]{siamltex}

\usepackage{amsmath}
\usepackage{amssymb}

\def\R{{\mathbb{R}}}

\def\R{{I\!\!R}}
\newcommand{\beas}{\begin{eqnarray*}}
\newcommand{\eeas}{\end{eqnarray*}}
\newcommand{\Be} {\begin{equation*}}
\newcommand{\Ee} {\end{equation*}}
\newcommand{\be} {\begin{equation}}
\newcommand{\ee} {\end{equation}}
\newcommand{\bea} {\begin{eqnarray}}
\newcommand{\eea} {\end{eqnarray}}
\newcommand{\pa} {\partial}

\newcommand{\ba} {\beta}

\newcommand{\Om} {\Omega}
\newcommand{\om} {\omega}

\newcommand{\la} {\lambda}

\newcommand{\noi} {\noindent}

\newcommand{\var} {\varepsilon}

\newtheorem{remark}[theorem]{Remark}

\title{ Unbounded Viscosity Solutions of Hybrid Control Systems}

\author{Guy Barles\thanks{Laboratoire Math\'ematique et Physique
Th\'eorique, F\'ed\'eration Denis Poisson, Universit\'e Fran\c{c}ois Rabelais Tours, Parc de Grandmont, 37200, Tours, {\sc FRANCE} ({\tt barles@lmpt.univ-tours.fr })}
\and Sheetal Dharmatti\thanks{Laboratoire MIP, UMR CNRS 5640, Universit{\'e} Paul Sabatier,
31062 Toulouse Cedex 9, {\sc FRANCE}({ email: sheetal@mip.ups-tlse.fr })}
\and Mythily Ramaswamy\thanks{TIFR Centre for Applicable Mathematics,
  Sharada Nagar, Yelahanka New Town, Bangalore-560065, {\sc INDIA} ({\tt mythily@math.tifrbng.res.in}).}}
\date{}
\begin{document}
\maketitle
\begin{abstract}
    We study a hybrid control system in which
both discrete and continuous controls are involved. The discrete
controls act on the system at a given set interface. The
state of the system is changed discontinuously when the
trajectory hits predefined sets, namely, an autonomous jump
set $A$ or a controlled jump set $C$ where controller can
choose to jump or not. At each jump, trajectory can move to
a different Euclidean space. 
We allow the cost functionals to be unbounded with certain growth
and hence the corresponding  value function
can  be unbounded.  We characterize the value function
as the unique viscosity solution of the associated quasivariational
inequality in a suitable function class.  We also consider the 
evolutionary, finite horizon hybrid control problem with similar model and 
prove that the  value function is the unique  viscosity solution in the 
continuous 
function class while allowing cost functionals as well as the dynamics to be 
unbounded.
\end{abstract}
\begin{keywords} 
Dynamic programming Principle, viscosity solution,
 quasivariational inequality, hybrid control 
\end{keywords}
\begin{AMS}  34H05, 34K35, 49L20, 49L25
\end{AMS}
\pagestyle{myheadings}
\thispagestyle{plain}
\section{ Introduction}
  The hybrid control systems are those involving
continuous and discrete dynamics and continuous and discrete 
controls. 
 Typical examples of such systems are constrained robotic
 systems \cite{bak} and automated highway systems \cite{var}.
For  some more 
examples of such systems, see
  \cite{bor},\cite{br} and references therein.

 In \cite{bor}, Branicky, Borkar and Mitter have 
presented  a model for a general hybrid control system in 
which continuous controls are present and in addition discrete 
controls act at a given set interface, which corresponds to 
the logical decision making process as in the above examples.
 The state of the system is changed discontinuously when the 
trajectory hits these predefined sets, namely, an autonomous 
jump set $A$ or a controlled jump set $C$ where controller can
 choose to jump or not. They  prove right  continuity of 
the value function corresponding to this hybrid control 
problem. Using dynamic programming principle they arrive at 
the partial differential equation satisfied by the
value function which turns out to be quasivariational 
inequality, referred hereafter as QVI.

 
  In an earlier work  \cite{SDMR1}, the authors have studied  this problem 
and  proved the local H\"older 
continuity of the value function and have characterized  it as
the unique viscosity solution of the QVI. There the cost functionals were 
assumed 
to be bounded and hence the value function was also bounded. This was 
essential in the
uniqueness proof, as the auxiliary function had to be modified a 
finite number of times to get the comparison result.
 In this paper we allow the cost functions to be unbounded. We use a different
 method to compare unbounded value functions, namely a suitable change of 
variable  to  reduce the unbounded function to a 
bounded one 
 (See  \cite{ba94}, for example ).
 We also tackle the uniqueness for the time
dependent problem, using a test-function  similar to the one introduced in
\cite{ba03}.

The paper is organized as follows:
In section 2, we introduce the notations, assumptions and quasivariational
inequality (QVI) satisfied by the value function. In section 3,
 we show that the value function is continuous. The next section  
  deals with uniqueness of the solution
 of QVI. We give a comparison principle proof  characterizing the
value function as unique viscosity solution of QVI ,  for the
stationary case in section 4 and then for the time dependent case in section 5.

\section{ Preliminaries }\label{prelim}

 In a hybrid control system, as in \cite{bor}, the state 
vector during 
continuous evolution,
is given by the solution of the 
following problem  
\begin{eqnarray} \label{1}
\dot{{X}} (t) &=& f_i ({X} (t), u (t))   \\
      {X} (0) &=& x \label{2}
\end{eqnarray} 
\noi where $ x \in  \Om_i $, a connected 
subset of $\R^{d_i},  d_i \in Z_+ $ for $i \in Z_1 \subset  Z_+  $.
Here  $f_i: \Om_i \times {\cal{U}}  \rightarrow  
\Om_i$  and  the continuous control set is :
 $$ \mbox {$\cal{U}$} =  
\left\{ u : [~0,\infty) \rightarrow U ~|~ u~\mbox{ measurable},\; U \; \mbox 
{compact metric space. }\;
 \right\} $$
 The trajectory  undergoes discrete jump  when it hits
 predefined sets $A_i $  the autonomous jump set and $C_i $, 
 the controlled jump set, both 
subsets of  $\R^{d_i}$.  The trajectory starting from $ x \in \Om_i$, 
on hitting
 $ A_i $, can jump to a  predefined 
destination set $D_j $  in another Euclidean space $\Om_j$ and continue the 
evolution there. This jump is given by
prescribed transition  map $ g:  A_i \times {\cal {V}}  \rightarrow  
\bigcup\limits_{i} D_i$, where $  {\cal {V}}$
is the discrete control set.  On hitting $C_i$ the controller can choose 
either 
to jump or not to jump. If the controller chooses to jump, 
then the trajectory is moved to a new point in
$D_k$, possibly in another space $\Om_k $. This gives rise to a sequence of
 hitting times of $\, \bigcup\limits_{j}A_j$,
 which we denote by $\sigma_i$ and sequence of hitting times 
of $\, \bigcup\limits_{j}C_j$, where the controller chooses to jump 
which is 
denoted by $ \xi_i $.

  We introduce the state space  
$$ \Om := \bigcup\limits_{i} \Om_i \times \{i\}, \quad i \in Z_1 \subset Z_+ ,
$$
and  the dynamics  $f: \Om \times \cal{U}  \rightarrow $ 
$\Omega$. Actually, $f= {f_i} $ with the understanding, 
$ \dot{{X}} (t) = f_i ({X} (t), u (t)) $ whenever $x \in \Omega_i$.
The predefined sets are 
$$
\begin{array}{lllll} 
A &=& \bigcup\limits_{i} A_i \times \{i\}
& A_i \subseteq \overline{\Om_i } \subseteq \R^{d_i}~; \\
C &=& \bigcup\limits_{i} C_i \times \{ i\}
& C_i \subseteq \Om_i \subseteq \R^{d_i}~; \\
D &=& \bigcup\limits_{i} D_i \times \{ i\}
&D_i \subseteq \Om_i \subseteq \R^{d_i}~. 
\end{array} 
$$ 

 The
trajectory of this problem, is composed of continuous evolution
 given by (\ref {1}) between two hitting times and discrete
 jumps at the hitting times. 
We denote $ X (\sigma_i^-, u(\cdot))$ 
by $x_i$, the point before an autonomous jump and $g(X(\sigma_i^-), v)$ by 
$x_i^\prime$, after the jump. The controlled jump destination of 
$ X(\xi_i^-, u(\cdot)) $  is
${X(\xi_i^-)}^\prime$, or simply, ${X(\xi_i)}^\prime$.
 For example, the dynamics for 
$\sigma_i < \xi_k  < \sigma_{i+1}  $, is given by
\begin{eqnarray*} 
\dot{{X}} (t) &=& f ({X} (t), u (t))   \quad \sigma_i < t < \xi_k , \\
      {X} (\sigma_i) &=&  {X} (\sigma_i^+) = g(X(\sigma_i^-), v)=
      x_i^\prime ,
\end{eqnarray*} 
and then 
\begin{eqnarray*} 
\dot{{X}} (t) &=& f ({X} (t), u (t))   \quad  \xi_k < t < \sigma_{i+1} , \\
      {X} (\xi_k) &=& X(\xi_k^+)= {X(\xi_k^-)}^\prime ,
\end{eqnarray*}

  We give the inductive limit topology on $\Om$ namely, 
$$ (x_n, i_n)  \in \Om \mbox{ converges to } (x,i) \in \Om  ~~
\mbox{ if  for some } N ~~\mbox {large and } \forall\; n \geq N $$
$$ i_n = i ~~ x, x_n \in \Om_i,\quad \Om_i \subseteq \R^{d_i}, \mbox{for some i,  and} ~\|x_n -x\|_{\R^{d_i}} < \var.$$
With the understanding of above topology we suppress the second variable $i$ from $\Om$. We follow the same for $ A, C$ and $D$. We also denote by $|x|$ the quantity which is equal to $\|x\|_{\R^{d_i}}$ if $x \in \overline{\Om_i}$.

We make the following basic assumptions on the sets 
$ A, C, D$ and on functions $ f \; \mbox{and}\; g $: \\
 \noi {\bf(A1)}: Each $\overline{\Om_i}$ is  closure of a connected, open subset of $\R^{d_i}.$ \\
\noi{\bf (A2)}: $A_i, C_i, D_i$ are closed, and for all $i$ and for
all 
$ x \in D_i$, $|x| < R$.
$\pa A_i $,$\pa C_i $ are $C^2$  
{and} $\pa A_i \supseteq  \pa \Om_i \quad \forall \;i. $  

\noi {\bf(A3)}: $g : A \times {\cal{V}} \rightarrow D$ 
is bounded,
 uniformly Lipschitz continuous map, with Lipschitz 
constant 
$G$ with the understanding that $ g = \{g_i\} $ 
and $g_i : A_i \times  {\cal{V}} \rightarrow \cup_j D_j $. 

\noi {\bf(A4)}: Vector field $f$ is Lipschitz continuous with 
Lipschitz constant $L$ in the state variable  $x$ and uniformly 
continuous in 
control variable $u$. Also, 
\be \label{3}
|f(x,u)| \leq F  \quad \forall \; x \in \Om ~\mbox{and}~~ 
\forall \; u \in U.
\ee

\noi {\bf(A5)}: Each ~$\partial A_i$ is compact for all $i$,
 and for some $\xi_0 > 0$, following transversality 
condition holds:
 \begin{equation}\label{9i}
 f(x_0,u)\cdot \zeta(x_0) ~\leq - 2 \xi_0 ~\;\forall~ \;x_0 ~\in ~\partial A_i ~~ 
\forall ~ u \in U
\end{equation} 
where $\zeta(x_0)$ is the unit outward normal to 
$\partial A_i$ at $x_0$.
We  assume similar transversality condition on $ \pa C_i$.

\noi {\bf(A6)}:  We assume that
\be \label{4}
\inf\limits_i d (A_i, C_i) \geq \beta 
\quad \mbox{and} \quad 
\inf\limits_{i} \; d (A_i, D_i) \geq \ba > 0 
\ee
 where $d$ is the appropriate Euclidean distance. Note that, 
above rules out infinitely many jumps in finite time.\\

\noi{\bf(A7)}: The control sets $ U $ and $ \cal{V} $  are assumed 
to be compact metric spaces.  

\medskip
  Now $ \theta : = (u(\cdot), v, \xi_i, {X (\xi_i)}^\prime)$ 
is the control and total discounted cost is given by 
\begin{eqnarray}\label{5}
J (x, u(\cdot), v, \xi_i, {X(\xi_i)}^\prime)
 &=&\int\limits_{0}^\infty K (X (t),u (t)) e^{-\lambda t}dt 
+ \sum\limits_{i=0}^\infty C_a (X (\sigma_i^-) , v )
 e^{-\lambda \sigma_i} \nonumber \\
 &+& \sum\limits_{i=0}^\infty C_c (X(\xi_i^-), {X(\xi_i)}^\prime)
 e^{- \lambda \xi_i} 
\end {eqnarray}
\noi where $\lambda$ is the discount factor, 
\noi $K : \Om \times \cal{U} \rightarrow $$\R_+$ is the running 
cost, $C_a : A \times {\cal{V}} \rightarrow$ $\R_+$ is the 
autonomous jump cost and $C_c : C \times D \rightarrow \R_+$ 
is the controlled jump cost. 
 The value function $V$ is then defined as: 
\begin{equation}\label{6}
V (x) = \inf\limits_{\theta \in 
({\cal{U}} \times {\cal{V}} \times [0, \infty) \times D)} 
J (x, u(\cdot), v, \xi_i, {X(\xi_i)}^\prime)
\end{equation}

  We assume the following conditions on the cost
 functionals:

\noi{\bf (C1)}: $K$ is nonnegative, continuous in the $x$ variable with at most
polynomial growth of degree $k$ , with $k$ and the growth being
independent of $i$.  $K$ is 
uniformly continuous in $u$ variable. 

\noi{\bf (C2)}: $C_a(x,v)$ and $C_c(x,x^\prime)$ are continuous in both
variables, uniformly continuous in $x$, uniformly with respect to $v$
and 
$x^\prime$ respectively, and bounded below by ${C^\prime} ~>~ 0$. 
Moreover $ C_a $ and
$ C_c $ have at most polynomial growth of degree $k$ in the first variable,
with $k$ and the growth being independent of $i$.
 
Note that under (C1) and (C2), value function is always non-negative
 and
 hence bounded below by 0. 
Using dynamic programming principle, one can show that the
value function satisfies the following QVI in viscosity sense:

\begin{theorem}\label{th:qvi} [{\em Quasivariational Inequality}]
Under the assumptions $(A1- A7)$, $(C1), (C2) $ and if $\lambda > kL$, the value function $V$ described in (\ref{6}) is continuous and has at most a polynomial growth of degree $k$. Moreover it satisfies the
following quasivariational inequality in the 
viscosity sense:
\begin{equation*}
V (x) = 
\begin{cases} MV (x) &\quad \forall \; x \in A  \\ 
 \min \left\{ NV(x), - H(x, DV(x)) \right\} &\quad \forall\; x \in  C \\ 
 -H (x, DV(x)) &\quad \forall\; x \in \Omega \setminus {A \cup C} 
 \end{cases} \eqno{(QVI)}
\end{equation*}
\\
\noi For $\phi,$ a function defined on $\Om$  and bounded below,
$M, N $ and $H$ are given by 
\beas
 M \phi (x) &=& \inf_{v \in {\cal{V}}} \{ \phi (g(x, v)) + C_a(x, v) \}, \\
  N \phi (x) &=& \inf_{x^\prime \in D } \{ \phi (x^\prime) + C_c(x, x^\prime) \}, \\
 H (x,p) &=& \sup_ { u \in U } \left\{ \frac { - K ( x, u) - f(x, u) \cdot p }{\la} \right\}.
 \eeas
\end{theorem}

The condition $\lambda > kL$ in Theorem~\ref{th:qvi} is used (and needed in the general case) to ensure that the value function is locally bounded and has a polynomial growth. For the uniqueness result for the quasivariational inequality (cf. Section~\ref{sec:uihc}), such condition is not playing any role.

Now we consider the hybrid control problem, with autonomous and
controlled jumps as before, but in finite horizon, namely for $t \in [0,T] $. In this case the dynamics depends on the time variable $t$ and we also allow cost functionals to depend on $t$.  We relax the global Lipschitzness assumption on $f$ and assume only local Lipschitzness.
Further relaxing boundedness of $f$, we allow it to grow linearly.

In order to simplify we assume that the cost $C_a$, $C_c$ are as above independent of time. It is not difficult to check that, besides additional technical details, the case when $C_a$, $C_c$ depend on $t$ can be treated in a similar way, with suitable adaptations.

The total cost in the time dependent case is given by,
\begin{eqnarray}\label{timedepeqn5}
J (s, x, u(\cdot), v, \xi_i, {X(\xi_i)}^\prime)
 &=&\int\limits_{s}^T K (t, X (t),u (t)) dt 
+ \sum\limits_{s \leq \sigma_i < T } C_a (X (\sigma_i^-) , v )
  \nonumber \\
 &+& \sum\limits_{s \leq \xi_i < T } C_c (X(\xi_i^-), {X(\xi_i)}^\prime) + h( X_x(T))
\end {eqnarray}
where $K : [0, T ]\times \Om \times \cal{U} \rightarrow $$\R_+$ is the running 
cost, $C_a :  A \times {{\cal V}} \rightarrow$ $\R_+$ is the 
autonomous jump cost, $C_c : C \times D \rightarrow \R_+$ 
is the controlled jump cost and $h$ is the terminal cost. 
 The value function $V$ is then defined as: 
\begin{equation}\label{timedepeqn1}
V (s, x) = \inf\limits_{\theta \in 
({\cal{U}} \times {\cal {V}} \times [s, T) \times D)} 
J (s, x, u(\cdot), v, \xi_i, {X(\xi_i)}^\prime)
\end{equation}

The precise assumptions are as follows: as before we assume  (A1)-(A3), (A6) and (A7). 
Because  $ f $ and $ K $ are time dependent we further assume the following:

\noi {\bf(A4)'}: Vector field $f(t,x,u)$ is locally Lipschitz continuous
 in the state variable  $x$ and uniformly 
continuous in the rest of the variables.
 Its linear growth is given by 
$$
|f(t, x,u)| \leq F (1 + |x|) \quad \forall \; x \in \Om ~\mbox{and}~~ 
\forall \; (t,u)  \in  [0,T] \times U.
$$

\noi {\bf(A5)'}: Each ~$\partial A_i$ is compact for all $i$,
 and for some $\xi_0 > 0$, following transversality 
condition holds:
 \begin{equation}\label{9ii}
 f(t, x_0,u)\cdot \zeta(x_0) ~\leq - 2 \xi_0 ~\;\forall~ \;x_0 ~\in ~\partial A_i ~~ 
\forall ~ (t,u )\in  [0,T] \times U
\end{equation} 
where $\zeta(x_0)$ is the unit outward normal to 
$\partial A_i$ at $x_0$.
We  assume similar transversality condition on $ \pa C_i$.
The cost functionals are assumed to satisfy 

\noi{\bf (C1)'}: $K$ is nonnegative and continuous in all the three variables;  moreover $K$ is uniformly bounded for bounded $|x|$.

One can derive the DPP for time dependent problem which is given as follows: 
 
For $s^\prime$ such that $ s \leq s^\prime < T$,
\beas
V (s, x) = \inf\limits_{(u, \xi_1)} \left[1_{s^\prime<  ( \sigma_1 \wedge \xi_1) } 
\left\{ \int\limits_{s}^{s^\prime} K (t, X (t),u (t)) dt + V ( s^\prime, X (s^\prime)) \right\} \right.  \\ 
 \left. + 1_{  \sigma_1 < (s^\prime \wedge \xi_1) } 
\left\{ \int\limits_{s}^{ \sigma_1 } K (t, X (t),u (t)) dt + V ( \sigma_1, X ({\sigma_1}^-)) \right\}  \right.  \\ 
 \left.  + 1_{  \xi_1 <  (s^\prime \wedge \sigma_1) } 
\left\{ \int\limits_{s}^{ \xi_1} K (t, X (t),u (t)) dt + V ( \xi_1, X ({\xi_1}^-)) \right\} \right] 
\eeas

It can be verified that the quasivariational inequality  (QVI-T) satisfied by $ V(s, x) $, is given by
$$
\left.
\begin{array}{rcl}
 V  - MV  &=& 0 \quad \hbox{in  }(0, T) \times A  \\ 
  \max \{ V - NV, -V_s  + H(s, x, DV) \} &=& 0 
\quad \hbox{in  } (0, T) \times C  \\ 
 -V_s + H (s, x, DV) &=& 0   \quad \hbox{in  } (0, T) \times \Omega \setminus {A \cup C}  \end{array}
 \right\}
 \eqno{(QVI - T)}
$$

\noi where for $ \phi$, a function defined on $([0,T] \times \Om)$
and bounded below, $M, N $ and $ H$ are given by 
\beas
M \phi(t, x) &=&  \inf_{v \in V } \{  \phi(t, g (x, v)) +  C_a (x, v) \} \quad \forall~ x ~\in A \\
  N \phi (t, x) &=&  \inf_{x^\prime \in D}\{ \phi(t,  x^\prime) + C _c (x, x^\prime)  \} \quad \forall ~ x ~\in C \\
 H ( t, x, p) &=& \sup_{ u \in U } \{ - K ( t, x, u) - f(t, x, u) \cdot p \} 
\eeas
For the terminal condition, it turns out to be more complicated than it is usually the case because of the possible (or imposed) jumps. In fact, the terminal data is obtained by solving the stationary problem
\begin{equation}\label{eqn-td}
\left.
\begin{array}{rcl}
V (T, x) - MV (T, x) &=& 0 \quad \forall \; x \in A  \\ 
\max \{ V(T, x) - NV(T, x), V (T, x) - h (x) \} &=& 0 
\quad \forall \; x \in C  \\ 
V( T, x)- h (x) &=& 0   \quad \forall \; x \in \Omega \setminus {A \cup C}.
\end{array}
\right\}
\end{equation}
 In order to have a continuous solution $V$ (and therefore a continuous terminal $V (T, x)$), we have to assume 

\noi{\bf (D1)}: Equation (\ref{eqn-td}) has a (unique) bounded from below solution $\tilde  h$ which is uniformly continuous for bounded $|x|$.

We point out that the key property in this assumption is the existence
of such 
uniformly continuous solution (uniformity means really uniformity in
$i$) 
while uniqueness is rather easy to obtain. This is really an
assumption on $h$.
 To go further in this direction, we also remark that if $C_c$ satisfies
$$C_c (x,y) \leq C_c (x,z) + C_c (z,y) \quad \hbox{for any  }x \in C, y \in D \hbox{ and }z \in C \cap D\; ,$$
which means that to jump once is always better that to jump twice, then $\tilde h$ can be built in the following way\\
(i) $\tilde{h}$ is given except on $A$ and $C$.\\
(ii) For points in $C$, set $\tilde{h} = min (h, Nh)$.\\
 The above condition implies that, actually, $Nh \leq N(Nh)$ on $C\cap D$ and therefore the right inequality holds in (\ref{eqn-td}) \\
(iii) For points in $A$, set $\tilde{h} = M\tilde{h}$. \\
This is well-defined since $\tilde{h}$ is known on $D$.

Assumption (D1) consists in assuming that this function $\tilde{h}$ is
uniformly continuous for bounded $|x|$, which is both a continuity
assumption on $h$ and a compatibility condition on its value near
$\partial A$ and $\partial C$. If for example,
$ C \cap D $ is empty,  then assuming $ h=Mh$ on $\partial A$
 and $h \leq Nh$ on
$\partial C$ in addition to the uniform continuity of $h$ on bounded subsets of $\overline{\Omega } $, will ensure that of  $\tilde{h}$, since $\tilde{h}=h$. One can construct easy examples in one dimension to show that if $ C \cap D $ is nonempty, the situation gets complicated and the conditions on $h$ are not
that transparent.

Once $\tilde h$ is known, we are left with a more classical terminal condition
 $$V( T, x)-\tilde  h (x) = 0   \quad \forall \; x \in \Omega.$$

\section{Continuity of the value function} 
Let the trajectory given by solution of $(\ref{1})$ and starting 
from the point $ x $ be denoted by $X_x(t,u(\cdot))$. Since 
$ x \in \Om$, in 
particular it
 belongs to some $ \Om_i$. 
Then we recall from theory of ordinary differential equations: 
\bea \label{7}
|X_x(t,u(\cdot)) -  X_z(t,u(\cdot)) | \leq e^{Lt}|x-z|, \\
|X_x(t,u(\cdot)) - X_x(\bar{t},u(\cdot)) 
\leq F |t - \bar {t}|,  \label{8}
|X_x(t,u(\cdot))|  \leq  |x| e^{Lt} + \frac{C}{L} ( e^{Lt} - 1)
\eea

For the derivation of these, see for example, (\cite{bardi} , Appendix, 
Chapter 5).
Here $F$ and $L$ are as in (A4). Now we can proceed as in 
\cite{SDMR1} using the transversality condition (A5) and conclude the 
continuity of the value function.
 Moreover, as running cost $K$, and discrete jump costs $C_a$ and $C_c$
 have polynomial growth of order $k$, one can see that total cost $J$
also has at most the same growth because of the assumption
$k L < \lambda $. Then it follows that the
  value function $V$ also has  a polynomial growth of order $k$,
  uniformly in all $\Omega_i  $.

For  $\eta > 0$, we introduce the following spaces of functions~: we
denote by 
$ E_\eta (\Om)$
the space of semicontinuous functions $u $ defined on $\Om$ 
and bounded below and such that, in each $\Omega_i $
$$ \lim_ {|x | \rightarrow \infty}u(x) e^{ - \eta |x|} = 0 ,$$
the limit being uniform w.r.t. $i$. Then $E(\Om)$ is defined by
  $$ E(\Om) = \cup_{\eta < \la /F} E_\eta \; . $$

In section 4, we prove the uniqueness of viscosity solution of QVI in
this function class $E$ to which the value function belongs. 

 In the case of finite horizon problems, the continuity of the value function is a more delicate problem : indeed, in the infinite horizon case, the fact that time is unbounded allows to use the transversality condition (A5) and to compensate the jumps, perhaps with a delay, and a priori this is not possible anymore in the finite horizon case. This shows that the problem is even more important when $t$ is close to the terminal time $T$ and this explains the above discussion about (D1) and the terminal condition $\tilde{h}$. 

To prove directly (and without to many technicalities) that the value
function is continuous requires (D1) but also the uniform continuity
on bounded subsets of $\Omega$ of the functions $K$, $C_a$, $C_c$,
which are rather strong 
assumptions. An alternative proof (but essentially with the same assumptions) is to show that the upper and lower semi-continuous envelopes of $V$ are respectively viscosity sub and supersolutions of (QVI-T) and to use the comparizon result of Section~5; this provides a slightly simpler proof.


 \section{Uniqueness - Infinite Horizon Problem}\label{sec:uihc}
In this section, we state and prove a comparison between sub and supersolutions of QVI, which belong to $E(\Om)$. 
\begin{theorem}
Assume $(A1)-(A7)$ and $(C1)-(C2)$. Let $u_1 , u_2 \in E
(\Om)$,  be respectively upper semicontinuous  subsolution and lower semicontinuous supersolution of the quasivariational inequality
given by (QVI) in the viscosity sense. Then $u_1 \leq u_2$ in $\Omega$.
\end{theorem}
\begin{proof}
We complete the proof  in 4 steps. In  step 1, we convert the 
unbounded value functions $ u_1, u_2$ 
into bounded functions $ w_1, w_2$ by a suitable 
change of variable.  We derive the 
modified quasivariational inequality  satisfied by $w_1$ and $w_2$. 
In the following 2 steps, for fixed $ \mu \in (0, 1)$, close to $1$,
we argue on $ \mu w_1$ and $w_2$ :
in Step 2, we examine the possibility of having approximate suprema at points of $ A \cup C$ and the actual consequences of such facts, while, in Step 3, we use the previous results to reach a contradiction.
In Step 4, we conclude the comparison of $ u_1 $ and $u_2$ using steps 2 and 3. 

\noindent{\bf{ Step 1}}:  Let $u_1, u_2$ be respectively upper and lower semicontinuous 
sub and supersolutions of the QVI in $E(\Omega)$.  
Fix $0< \eta < \frac{\la}{F}$ such that both $u_1, u_2 $ lie in $E_{\eta}$, where, as above, $\la$ is the discount factor
and $F$ is the bound on dynamics $f$. 
Define $$ w_1 (x) = u_1 (x) e^{ - \eta \xi(x)}, \quad  w_2 (x) = u_2 (x)  e^{ - \eta \xi(x)},  $$
where $ \xi: \Om \rightarrow R^+ $ is  a smooth function such that,  
$ | D \xi | \leq 1  $ for all $x \in \Om$,  
 $ \xi (x ) = 0 $ if  $|x| \leq R $ and $ \xi (x) = (1 + |x|^2)^{1/2}$ if $ |x| > 2R  $ where $R$ is given by (A2).
Since $\xi$ behaves like $|x|$ at infinity and $u_1, u_2  \in E_{\eta}$, $ w_1 (x), w_2 (x)$ are bounded and converge to $0$ when $x$ tends to infinity in each $\Omega_i$, uniformly w.r.t. $i$.

The functions $w_1, w_2$ are respectively upper and lower semicontinuous  sub and supersolution of a QVI which can be first written as
\begin{equation*}
w (x) e^{ \eta  \xi(x)} = 
\begin {cases} M ( w (x) e^{ \eta\xi(x)}  )  &\quad \forall \; x \in A \\ 
 \min \left\{N ( w (x) e^{ \eta \xi(x)}), - H(x, D (w (x) e^{ \eta \xi(x)})) \right\} &\quad \forall\; x \in  C \\  
-H (x, D(w (x) e^{ \eta \xi(x)})) &\quad \forall\; x \in \Omega \setminus {A \cup C} 
\end{cases} 
\end{equation*}
But, recalling the facts that $|x| < R$ if $x \in D$ and that $g$ maps $A \times {\cal{V}}$ into $D$, we can simplify this QVI in the following way
\begin{equation*}
w (x) =
\begin {cases}  
 \tilde M w (x) &\quad \forall \; x \in A \\ 
\min \left\{   \tilde N w (x) , -\tilde H(x, w(x), D w (x) ) \right\} &\quad \forall\; x \in  C \\
-\tilde H (x, w(x), Dw (x) ) &\quad \forall\; x \in \Omega \setminus {A \cup C} 
\end{cases}
\end{equation*}
where
\begin{eqnarray*}
 \tilde M w (x) &=&
e^{- \eta \xi(x)}  \inf_{v \in \cal{V} } \{  w(g (x, v)) +  C_a (x, v) \}  \\
 \tilde N w (x) &=& e^{ - \eta \xi(x)} \inf_{x^\prime \in D}\{ w( x^\prime) + C _c (x, x^\prime)  \}
\\  
\tilde H(x, w, D w (x) ) &=& \frac{1} {\la}\left\{ 
\sup_{u\in U} \{ - f(x, u) \cdot ( Dw (x) + \eta  w(x) D \xi(x) ) \} 
- e^{ -\eta \xi(x)} K(x, u) \} \right\}
\end{eqnarray*}
Indeed, by definition of $\xi$ and for $i=1,2$, $w_i (g (x, v))= u_i (g (x, v))$ for any $x\in A$ and $v \in \cal{V}$, and, in the same way, $w_i (x^\prime)= u_i (x^\prime)$ for any $x^\prime \in D$.

Using the definition 
of $\tilde H $ and the assumptions on the cost functionals and dynamics we can show that 
\bea  \label{unbdeqn5}
 \nonumber  | \tilde H( y, w_2(y ), p_2 )- \tilde H( x , w_1(x ), p_1) | &\leq& 
 ( L/ \la) | x - y| | p_1|  + (F/ \la) | p_1(x) - p_2(y)|  \\ &+&
  (\eta F/\la ) | w_1(x) - w_2 (y)|   + (1/\la) \om_{\tilde K} ( |x - y|) \\ \nonumber
  &+&  (\eta/\la) | w_2| \sup_{u\in U} \{| D \xi(x)\cdot f(x, u) - D
 \xi(y)\cdot f(y, u)| \}
 \eea
where $ \tilde K (x, u) = e^{- \eta \xi (x)} K(x, u)$ and $ \om_{\tilde K} $ is the modulus of continuity of $\tilde K $.

Now we fix $ \mu \in (0, 1)$, close to $1$, and prove the comparison result for $ \mu w_1 $ and $w_2$. To do so, we assume, by contradiction, that $ \sup_{\bar \Om} (\mu w_1(x) - w_2(x)) > 0$ or, in other words
$$\sup_{j}\sup_{{\Om_j} } (\mu w_1(x) - w_2(x)) = m > 0 .$$ 
We first recall that $w_1 (x), w_2(x)$ converge to $0$ when $x$ tends to infinity in all $\Omega_j$, uniformly w.r.t. $j$; therefore each  $\sup_{{\Om_j} }$ is indeed achieved at some point. But since the set of indices $j$ is infinite, we have anyway to argue with approximate supremums (of course, if this supremum is attained at some finite point then the following proof gets simplified). We consider $\kappa>0$ small enough (a more precise estimate on its size will be given later on) and let $i$ and $x_{\kappa}$  be such that 
\begin{equation}\label{unbdeqn1}
 \sup_{\Om_{i}} (\mu w_1(x) - w_2(x)) = \mu w_1(x_\kappa) - w_2(x_\kappa)
 \geq  m - {\kappa}>  \frac{m}{2} > 0 .
\end{equation}
Note that $x_\kappa$ remains bounded : indeed we can find $R_m$ such that 
$$ |\mu w_1(x)| < m/10  \; \text{ and } \; |w_2(x)|  < m/10 \; \text{ for} \; |x| > R_m  \; \hbox{for all $x \in \Om_j$ and for all $j$.} $$
Then, $x_\kappa$ being the maximum of  $\mu w_1 - w_2$  in $\Om_i$ and
  $  m/2  <  m - \kappa < (\mu w_1 - w_2)(x_\kappa)$ will imply that 
$|x_k| < R_m.$

Next we are going to consider three different cases, namely, 
\begin{enumerate}
\item  $ x_\kappa \in A_{i} $,  
 \item  $ x_\kappa \in C_{i}$ and $w_2 (x_\kappa) \geq \tilde N w_2 (x_\kappa)$
 \item  Either $ x_\kappa \in \Om_{i} \setminus (A_{i} \cup C_{i})$ or $ x_\kappa \in C_{i}$ and $w_2 (x_\kappa) < \tilde N w_2 (x_\kappa)$.
\end{enumerate}
We consider cases 1 and 2 in the next
 step and case 3 in Step 3. \\

\noindent
{\bf{ Step 2:}} We first consider the case 1 namely $ x_\kappa \in
A_{i} $. 
For $w_1, w_2$ , viscosity  sub  and supersolutions, the 
conditions on $\pa A$  are  satisfied  only in the
viscosity sense: For the subsolution 
\be
\min \{w_1 +\tilde H (x, Dw_1, w_1), \quad w_1- \tilde M w_1\} \leq 0 
\quad \mbox{ on } \pa A \label{Eq1} \ee
 For the supersolution, 
\be
 \max \{w_2 +\tilde H (x, Dw_2, w_2), \quad w_2 -\tilde M w_2\} \geq 0 \quad
\mbox{ on }  \pa A. \label{Eq2} \ee
 Hence we need to rule out the bad inequalities for $w_1, w_2,$ namely 
$$w_1 > \tilde M w_1  \; \mbox{ or } \;  w_2 < \tilde M w_2 $$
on $\pa A $. Our assumption $(A5)$ helps us to avoid the above situation and 
to conclude that the right inequalities for $w_1, w_2$ hold on $\pa A.$ 
To do so, we need the following
\begin{lemma}\label{MandN}
 If $w_1$ and $w_2$ are viscosity sub and supersolutions 
of QVI and $w_1, w_2$ are u.s.c and l.s.c   on $\overline{\Om}$ respectively, 
then  the functions $\tilde M w_1, \tilde N w_1$ and $\tilde M w_2,
\tilde N w_2$ are respectively upper and lower semi-continuous. 
Moreover, for any $x\in \partial A$, we have 
$$w_1 (x) \leq \tilde M w_1(x) ,  \quad w_2 (x) \geq \tilde M w_2(x)$$ 
and for any $x \in \partial C$, we have $ w_1 (x) \leq \tilde N w_1 (x)$.
\end{lemma}

\medskip

We postpone the proof of this lemma to the end of the section.

We now conclude the proof by using Lemma~\ref{MandN}. If $x_\kappa \in
A$ 
then, using Lemma~\ref{MandN} for boundary points,  we have
\beas
w_1(x_\kappa)  \leq \tilde M w_1(x_\kappa) \quad \hbox{and}\quad 
w_2(x_\kappa) \geq  \tilde M w_2(x_\kappa)
\eeas
By definition of $ \tilde M$, and  compactness of control set ${\cal {V}}$, we can find $ v_0$ such that,
$$
\tilde M w_2(x_\kappa) = e^{ - \eta \xi (x_\kappa)}  \{  w_2(g (x_\kappa, v_0)) +  C_a (x_\kappa, v_0) \}
$$
For all controls $ v \in {\cal{V}}$ and hence in particular for $v_0$ we have,
$$
\tilde M w_1(x_\kappa) \leq e^{ - \eta \xi (x_\kappa)}   \{  w_1(g (x_\kappa, v_0)) + C_a (x_\kappa, v_0) \}
$$
Hence,
\beas
 \mu w_1(x_\kappa) - w_2(x_\kappa) & \leq & e^{ - \eta \xi(x_\kappa)} \{ \mu w_1(g (x_\kappa, v_0)) -   w_2(g (x_\kappa, v_0))  - ( 1- \mu ) C_a (x_\kappa, v_0)
\} \\
&\leq& e^{- \eta \xi (x_\kappa)} \{ \mu w_1(g (x_\kappa, v_0))- w_2(g (x_\kappa, v_0))- (1- \mu ) C^\prime\} \\
&\leq& e^{- \eta \xi (x_\kappa)} \{m -(1- \mu ) C^\prime \}
\eeas
This inequality and  (\ref{unbdeqn1})  imply 
$$ m - \kappa   \leq \mu w_1(x_\kappa) - w_2(x_\kappa) \leq e^{- \eta \xi (x_\kappa)} \{m -(1- \mu ) C^\prime \}
\leq m -(1- \mu ) C^\prime
$$
This is a contradiction if we choose $\kappa$  such that  
\be \label{eq1kappa}
( 1 - \mu)  C^\prime> \kappa.
\ee
Thus, for such a choice of $\kappa$, $ x_\kappa$ cannot be on $A_i$.

Next we consider the case when $ x_{\kappa} \in C_i $. If $w_2 (x_\kappa) \geq \tilde N w_2 (x_\kappa)$, we use Lemma~\ref{MandN} to deduce $w_1 (x_\kappa) \leq \tilde N w_1 (x_\kappa)$ and a similar argument as above leads to a contradiction under the same condition on $\kappa$. Therefore we may assume without loss of generality that $w_2 (x_\kappa) <\tilde N w_2 (x_\kappa)$ which means that we are in the third case.

\noindent {\bf{ Step 3:}} We examine here the third (and last possible) case. To do so, we introduce the function $\Phi :
\Om_i \times \Om_i \to \R $ by
\Be
\Phi (x, y) =   \mu w_1(x) - w_2(y) - \frac{1}{\epsilon}|x-y|^2 - |x-x_\kappa|^2
\Ee 
where $\epsilon$ is a small positive parameter, devoted to tend to $0$. Since $ \Phi  $ goes to $ - \infty$ as $|x|$ goes to $\infty  $, there exists a maximum point $(x_0, y_0)$ of $\Phi$.
Note that,
\begin{equation} \label{unbdeqn3}
\frac{m}{2}   <   \mu w_1(x_\kappa) - w_2(x_\kappa)  = \Phi( x_\kappa, x_\kappa) \leq \sup \Phi(x, y) 
= \Phi(x_0, y_0)\; .
\end{equation}
The following estimates, listed in a lemma,  are standard in the theory of viscosity solutions and 
can be derived from the definition of $ \Phi$ and $ (x_0, y_0)$.
For example see \cite {br}, Chapter 3, Theorem 2.12. Only the last property stated in the lemma, is non-standard. 
\begin{lemma}\label{lem5.1}
\begin{enumerate}
\item[(i)] ${|x_0-y_0|} \le\sqrt{C \epsilon} $ for some constant ${C} $ depending only on the $L^\infty$ bounds on  $w_1$ and $w_2$. 
\item[(ii)] $\displaystyle \frac{|x_0-y_0|^2}{\epsilon} \to 0$ as $\epsilon \to 0$
\item[(iii)] $x_0, y_0 \to x_\kappa$,  $w_1(x_0) \to w_1(x_\kappa)$, $w_2(y_0) \to w_2(x_\kappa)$ when $\epsilon \to 0$,
\item[(iv)] If $x_\kappa \in C_i$ and $w_2 (x_\kappa) < \tilde N w_2 (x_\kappa)$, $w_2 (y_0) <\tilde N w_2 (y_0)$ if $\epsilon$ is small enough.
\end{enumerate}
\end{lemma}

\medskip

We leave the proof of the point (iv) of this lemma to the reader: it just uses the fact that $w_2(y_0) \to w_2(x_\kappa)$ when $\epsilon \to 0$ together with the lower semi-continuity of $N w_2$ coming from Lemma~\ref{MandN}.

Now define the test functions $ \phi_1 $ and $\phi_2$ by
\Be
\phi_1(x) = w_2( y_0) +   \frac{1}{\epsilon}|x - y_0|^2 + |x-x_\kappa|^2 \quad \hbox{and}\quad
\phi_2(y) = \mu  w_1(x_0) - \frac{1}{\epsilon}|x_0 - y|^2 - |x_0 -x_\kappa|^2
\Ee
The function $\mu w_1 - \phi_1$ attains its maximum at $x_0$ and, for $\epsilon$ small enough, $x_0 \in \Omega_i \setminus A_i$ because $x_0 \to x_\kappa \in \Omega_i \setminus A_i$; therefore QVI implies
$$
\mu w_1(x_0) + \sup_{u\in U} \frac{1} {\la} \{ - f(x_0,u) \cdot (\mu D \phi_1 (x_0)
   + \eta \mu w_1 (x_0) D \xi (x_0) ) -\mu e^{- \eta \xi (x_0)} K (x_0,u) \}  \leq 0 \; .
$$
On the other hand, $ w_2 - \phi_2$ attains its minimum at $ y_0 \in \Omega_i \setminus A_i$ and $w_2 (y_0) <\tilde N w_2 (y_0)$ if $y_0 \in C_i$; therefore, by using definition of viscosity supersolution, we have
$$ w_2(y_0) + \frac{1} {\la}
\sup_{u\in U} \{ -  f(y_0, u)\cdot ( D \phi_2(y_0) + \eta w_2(y_0) D \xi(y_0) )  - e^{ -\eta
  \xi(y_0)} 
K(y_0, u) \}   \geq 0 \; .$$
Denoting by $\displaystyle p_\epsilon = \frac{2(x_0 - y_0) }{\epsilon} $ and substituting $ D \phi_1(x_0) =  p_\epsilon + 2(x_0 - x_\kappa)$, $D \phi_2(y_0 ) = p_\epsilon$ in the above and estimating as in (\ref{unbdeqn5}), we get, 
\begin{eqnarray*}
( \la - \eta F) ( \mu w_1(x_0) - w_2(y_0)) & \leq & L | x_0 - y_0| \{ | p_\epsilon|  + \eta | w_2| (F + L)| x_0 - y_0| \\
&&+ ( 1 - \mu ) \om_{\tilde K } (|x_0 - y_0|) + 2F |x_0 - x_\kappa| \; .
\end{eqnarray*}
Now remembering $ \la - \eta F  > 0, $
\beas
( \la - \eta F)\Phi( x_0 , y_0) &=&  ( \la - \eta F)\{ \mu w_1(x_0) - w_2(y_0) - \frac{1}{\epsilon}|x_0 - y_0|^2  \} \\
&\leq& ( \la - \eta F) ( \mu w_1(x_0) - w_2(y_0) ) \\
&\leq&  L | x_0 - y_0| \{| \frac {2(x_0 - y_0) }{\epsilon} \}
 + \eta ( F +  L) |w_2| | x_0 - y_0|  + \\
 & & ( 1 - \mu) \om_{\tilde K } (|x_0 - y_0|) + 2F |x_0 - x_\kappa| 
\eeas 
Using the lemma,  we can choose $ \epsilon  $  such that RHS of the above 
inequality can be made arbitrarily small, which will contradict 
 (\ref{unbdeqn3}) namely,
$$ \Phi( x_0 , y_0) \geq \frac{m}{2 } >0. $$

Thus, in each of the three steps, we arrive at a contradiction to (\ref{unbdeqn1}) and we conclude that,
$$
 \sup_i \sup_{\Om_{i}} (\mu w_1(x) - w_2(x))  \leq   0 .
$$
Now sending $\mu$ to 1 we get the required comparison between $ w_1 $ and $ w_2$, namely,
 $  w_1(x) \leq w_2(x) $ for all $x$ in $  {\Om}$. 
 
\noindent{\bf Step 4: } From Step 3, we have concluded that
$$  w_1(x) - w_2(x) \leq 0 \quad \forall ~ x \in \Om $$
Now, by definition of $w_1, w_2$,
$$ w_1 (x) = u_1 (x) e^{ - \eta \xi(x)}, \quad  w_2 (x) = u_2 (x)  e^{ - \eta \xi(x) },  $$
hence,
$$
 u_1(x) - u_2(x) =   e^{  \eta \xi(x)} \{ w_1(x) - w_2(x) \} \leq 0 \quad \forall ~ x \in \Om 
 $$
Thus, we have the comparison between $u_1$ and $ u_2, $ which are solutions
of the original quasivariational inequality.

In order to complete the proof, we have to provide the 

{\bf Proof of Lemma~\ref{MandN}}. 
This  uses the idea in lemma 5.2 (page 113) of \cite{ba94}.
 Suppose that for some $x \in \pa A_i$
$$ w_1 (x) >   \tilde{M} w_1 (x). $$
Let the distance function from $ A_i $ be denoted by $d(\cdot)$. 
Consider the test   function,  
$$ \phi (y) = w_1 (y) - \frac{|x -y|^2}{\var} - C d(y) $$ 
for positive parameters $\var $ and $C,$ to be fixed suitably later on. 
Let $y_\var$ be the max of $\phi$ in $\overline{\Om_i}.$ Then in 
particular,
 $$w_1 (y_\var) \geq w_1 (x)$$
for each $\var$ and $y_\var \rightarrow x. $ Hence 
$$\liminf w_1 (y_\var) \geq w_1 (x).$$ 
Using the uppersemicontinuous property of $w_1,$ 
$$ 
\limsup  w_1 (y_\var) \leq w_1 (x).$$ 
Thus 
$$\lim w_1 (y_\var) = w_1 (x)$$
By our assumption $w_1 (x) > \tilde{M} w_1 (x).$  Since $\tilde{M} w$ is 
u.s.c, we have 
$$ w_1 (y_\var) > \tilde{M} w_1 (y_\var)$$ 
for $\var $ sufficiently small. If $y_\var \in \pa A_i,$ then we have by 
(\ref{Eq1}), 
$$ 
w_1 (y_\var) + \sup\limits_{u} \frac{1} {\la} \{ -\langle f (y_\var, u), 
\frac{2 (x-y_\var)}{\var} - C n  (y_\var) 
+ \eta  w_1 (y_\var) D \xi(y_\var) \rangle  
- e^{ -\eta \xi(y_\var)} K (y_\var, u) \} \leq 0  $$
where for the distance function $d (\cdot)$ from $A,$ $ D d (x)= n (x). $  
Now, we divide by $C$ and let $\var \rightarrow 0$ and $C \rightarrow 
\infty,$ in  such a manner that $C  \sqrt{\var} \rightarrow 
\infty. $ Recalling that $ \frac{|x-y_\var|}{\sqrt{\var}} 
\rightarrow 0, $ we get 
$$ 
\sup\limits_u \{ \langle f(x, u), n(x) \rangle \} \leq 0,\quad \Rightarrow \quad
  \langle f (x, u), n (x) \rangle \leq 0 \quad \forall \quad  u \in {\cal{U}}. $$
As $ n (x) = - \zeta (x),$ appearing in $(A5)$ this conclusion contradicts 
$(A5).$ Hence our assumption cannot hold. 

\par In a similar manner arguing with  the test function. 
$$ 
\phi (y) = w_2 (y) + \frac{|x -y|^2}{\var} + cd (y)$$
we get the required inequality of $w_2.$ 
Thus the proof is complete.
\end{proof}

This  theorem characterizes  the value function of hybrid control problem 
as the unique viscosity solution of the QVI in the function class $E$ defined 
earlier. We would like to remark the following.

\begin{remark}
We can extend the uniqueness result to the hybrid game theory problem with unbounded value functions. 
 The case of bounded value functions was treated in \cite {SDMR2}.
 \end{remark}

\section{Uniqueness - Finite Horizon Problem}

The comparison proof for finite horizon problem is proved by building test-functions, which go to infinity at the boundary of the domain of comparison and tend to zero in a specific subset of the interior. These are similar to the ones in Ley \cite{OL}  and Barles, Biton and Ley\cite{ba03} where they are introduced as ``friendly giants'', either for proving ``finite speed
of propagation type results'' for first-order equations or for proving that the comparison holds for second-order equations. Here it is extended to quasi-variational inequalities involving HJB equations.

First we announce the uniqueness theorem.
\begin{theorem}-
Assume $(A1-A3), (A4)', (A5)', (A6), (A7)$ and $(C1)' , (C2)$. Let
$u_1 , u_2 $ be respectively a u.s.c. subsolution  bounded from below
and a l.s.c supersolution  bounded from below of the quasivariational 
inequality (QVI-T) in $(0,T) \times \Omega$. If $u_1, u_2$ satisfy: for any $r>0$, $u_1, u_2$ are bounded on $[0,T]Ê\times \{x \in \Om\,;\, |x| \leq r\}$ and
\be\label{hypid}
\limsup_{t\to T}\, \sup_{|x|\leq r}\,\left [ u_1 (t,x) - u_2 (t,x)\right ] \leq 0 \; ,
\ee
then $u_1\leq u_2$ in $(0,T) \times \Omega$.
\end{theorem}

\medskip

We point out that (\ref{hypid}) is a non trivial assumption since we work in $\Om = \bigcup\limits_{i} \Om_i \times \{i\}$ because of the dependence in $i$; roughly speaking, it means that the terminal data is assumed with a certain uniformity w.r.t $i$. This assumption is satisfied by the value function of the control problem under mild additional assumptions on the data as we show it at the end of the section.

\begin{proof}  Since $u_1 , u_2 $ are both bounded from below, we can assume as well that they are positive since we can add the same positive constant to both of them and the resulting functions still solve (QVI-T) in $(0,T) \times \Omega$ (of course, (\ref{hypid}) also holds.)

We want to prove that for fixed $\mu $ close to 1 in $(0,1)$,
$$w = \mu u_1 - u_2 \leq 0 \; \mbox{ in } \; [0,T]  \times \Omega \ .
$$
For that, we divide the interval $ [0, T] $ into finite number of sub-intervals
of  length less than a certain constant and prove the comparison on
this subinterval first and then extend the result to $ [0, T]$ by
repeating the argument on each of the subintervals. 

Let us set $\bar K := F(1+R+\beta)$, where $F$ is as in $(A4)'$, $R$
is the bound for $D$, introduced in $(A2)$ and $\beta $ is as in $(A6)$.
We choose $T_0$ such that
\be  \label{timechoiceeq1}
T- T_0 \, <  \, \min \{  \frac{1}{2F} , \frac{\beta}{2 \bar K }  \}
\ee
The first step consists in proving that $w \leq 0$ on $[T_0,T] \times
 D$.
To do so, we introduce a  neighbourhood of $D$, namely
$$  \tilde D := \{ x  \in \Om \ :\ d(x,D) < \beta \ \}$$
for $\beta >0$  as before and we set
$$ \tilde m(T_0) := \sup_{[T_0,T] \times D} w(t,x) .$$ 
Because of assumption (\ref{hypid})  and since 
$$w(t,x) = (\mu -1)u_1(t,x)+ (u_1 - u_2)(t,x) \leq (u_1 - u_2)(t,x)
,$$  
there exists  $\tau$, $\, T_0 \leq \tau \leq T $ such that
\be  \label{timechoiceeq2}
\tilde m ( \tau ) \,  \leq ( 1 - \mu ) \frac{C'}{2}   \,
\ee
where $ C' $ is the lower bound for the cost functionals introduced
in $\bf {(C2) } $.
We will show that $w$ is a subsolution of the variational inequality
\be \label{timedepeq1}
 \min \{ -w_t (t, x) - F ( 1 + |x|) |Dw | \ , \  w (t,x) \}
= 0 \;
 \ee
in 
$ \tilde D_{\tau} :=  (\tau, T ) \times \tilde D $ and from that we
will conclude 
 $ w \leq 0$ on $[ \tau, T ] \times D$
by using ``finite speed of propagation type properties'' 
(cf. for example Ley \cite{OL}). Once this is true for every such $\tau$,
it will follow that $w \leq 0$  in  $[T_0,T] \times D$.

From this first result, we will deduce, in a second step, that $w \leq
0$ on $[T_0,T] \times A$. Finally the
third step will be devoted to proving first that $w$ is a subsolution of
(\ref{timedepeq1}) in $[T_0,T] \times \Omega \setminus A$
and then $w \leq 0$ in $[T_0,T] \times  \Omega$.\\

\noindent{\bf Proof of Step 1 :}  Let $\Phi \in C^2 ( (0, \infty) \times \Om_\kappa)$ and let $( \bar t, \bar x) \in \tilde D_{\tau}$ be a strict local maximum point of $ w - \Phi =  \mu u_1 - u_2 - \Phi $ in $\tilde D_{\tau} $. We want to prove that
$$  \min \left\{  -\Phi _t (\bar t,\bar x) - F ( 1 + |\bar x|) |D\Phi (\bar t,\bar
x) | \ ,\ w( \bar t, \bar x) \right\} \leq  0 \; .$$
If $w( \bar t, \bar x)  \leq 0 $, there is nothing 
to prove. 

If $w( \bar t, \bar x) > 0 $, since $C \cap D $ need not be empty, we
 first remark that
\be \label{eq-for-u-Nu}
u_2 (\bar t,\bar x) < N u_2 (\bar t,\bar x)\; \mbox{ if } \; \bar x \in C_\kappa
\ee
Indeed, otherwise, since $u_1 (\bar t,\bar x) \leq  N u_1 (\bar t,\bar x)$, we would have, for some $x'\in D$
$$  (\mu u_1 - u_2)(\bar t,\bar x) \leq (\mu u_1 - u_2)(\bar t,x') + (\mu-1) C_c (\bar x,x') \leq \tilde m (\tau) +  (\mu-1) C',$$
which contradicts the fact that $w(\bar t, \bar x) > 0 $, because of
our choice of $\tau$, satisfying (\ref{timechoiceeq2}).

For $ \var > 0 $, we consider for some $\rho$ such that
$(\bar t - \rho , \bar t + \rho ) \times B( \bar x, \rho) 
\subset  \tilde D_{\tau}$,
$$ \max_ {(\bar t - \rho , \bar t + \rho ) \times B( \bar x, \rho)^2 } 
\{  \mu u_1 (t, x) - u_2(t, y) - \Phi (t, x) - 
\frac  {| x-y|^2}{\var^2}\} \ .
 $$
It can be shown that the above maximum is achieved at points $(x_\var,
y_\var, t_\var) $ such that 
 $$
 (t_\var, x_\var, y_\var) \rightarrow ( \bar t, \bar x, \bar x) \quad \mbox{and}\quad \frac 
 {| x_\var - y_\var |^2}{\var^2} \rightarrow 0 \; \mbox {as} \; \var \rightarrow 0 ,
 $$
$$ u_2 (t_\var, y_\var)  \rightarrow u_2 (\bar t,\bar x). 
$$
Then, thanks to (\ref{eq-for-u-Nu}) and the lower semicontinuity of $N
u_2 $ 
and the above
convergence, we have, for $\var$ small enough,
$$ u_2 (t_\var, y_\var)  < N u_2 (t_\var, y_\var).
$$
Therefore, at these points $( t_\var, y_\var)$, $u_2$ is a viscosity 
supersolution of the HJB equation 
given by
$$- v_t + \sup_{ u \in U } \{ - K ( t, x, u) - f(t, x, u) \cdot Dv
 \} = 0 \; ,$$ 
while $v = \mu u_1$ is a viscosity subsolution of the (slightly different) HJB equation
$$- v_t + \sup_{ u \in U } \{ - \mu K ( t, x, u) - f(t, x, u) \cdot Dv
 \} = 0 \; .$$ 
By the results of the User's guide \cite{CIL}, there exists $a \in \R$
such that
\beas
 (a,p_\var) &\in & D^- u_2 (t_\var, y_\var) \\
 (a + \Phi_t (t_\var, x_\var), p_\var + D\Phi (t_\var, x_\var)) &\in &
D^+ u_1 (t_\var, x_\var)
\eeas
where  $\displaystyle p_\var := \frac{2(x_\var-y_\var) }{\var^2}$ and therefore
\beas
 -a - \Phi_t (t_\var, x_\var) + \sup_{ u \in U } 
\{ - \mu K ( t_\var, x_\var, u) - f(t_\var, x_\var, u) \cdot 
 (D \Phi(t_\var,x_\var )+p_\var) \} & \leq & 0 \; ,   \\
- a + \sup_{ u \in U } 
\{ -  K ( t_\var, y_\var, u) - f(t_\var, y_\var, u) \cdot 
p_\var  \}   & \geq  & 0   \; .   
\eeas
Subtracting these two viscosity inequalities, we get
\beas
 -\Phi_t (t_\var, x_\var)& - & \sup_{ u \in U } 
\{ -  K ( t_\var, y_\var, u) - f(t_\var, y_\var, u) \cdot 
p_\var  \} \\
& + & \sup_{ u \in U } 
\{ - \mu K ( t_\var, x_\var, u) - f(t_\var, x_\var, u) \cdot 
(D \Phi(t_\var,x_\var )+p_\var)  \}  \leq 0 
\eeas
Using that $\sup (\cdots)-\sup (\cdots) \geq \inf (\cdots - \cdots)$
and the 
continuity properties of $K$ and $f$, we are lead to
$$- \Phi_t (t_\var, x_\var) -  F (|1 + |x_\var|) | D\phi(t_\var,x_\var
)|
- \mu \om_K (|x_\var - y_\var|)
- \om_f ( |x_\var - y_\var ||p_\var|)\leq 0
$$
Thus, as $ \var \rightarrow 0$, noticing that $|x_\var - y_\var ||p_\var|= 2\frac{| x_\var - y_\var |^2}{\var^2} \to 0$, we obtain
\beas
- \Phi_t (\bar t, \bar x) - F (1 + |\bar x|) | D\phi(\bar t, \bar x
)|   
\leq 0
\eeas
proving the claim that $ w =  \mu u_1 - u_2 $ is a viscosity subsolution of 
(\ref{timedepeq1}) in $ \tilde D_{\tau} $.

Next we deduce that $w \leq 0$ on  $[\tau, T]Ê\times D$.
We first remark that, by {\bf (A2)}, if $x \in  \tilde D$, $|x| \leq R
+\beta$. Recall that
 $\bar K := F(1+R+\beta)$. Hence  $w$ is as well a subsolution of
\be \label{timedepeq2}
 \min  \{-w_t (t, x) - \bar K |Dw | ,  w(t,x) \} = 0 \quad\hbox{in }
\tilde D_{\tau} \; .
\ee
We now construct a strict supersolution for this equation
using ``friendly giant'' functions. For that, we consider an
increasing 
 $C^\infty$  function $\psi : (-\infty,\beta)
\to \R$ such that 
$$ 
\psi(s) \equiv 0 \; \mbox{ if} \; s\leq \beta/2 \quad \mbox{ and} \quad
 \psi(s) \to + \infty  \; {\mbox as} \; s \to \beta.
$$
 We claim that, for any $\kappa$, $x_0 \in D_\kappa$, $\eta >0$ and 
$T_0$ satisfying (\ref{timechoiceeq1}),
the function 
$$\chi_{x_0} (t,x) := \psi( \bar K (T-t) + |x-x_0|) + \eta (T-t)$$ 
is a strict (smooth) supersolution of (\ref{timedepeq2}) in the domain 
$${\mathcal C}_{x_0} := \{ (t,x) \in \tilde D_{\tau} \ : \ \bar K (T-t)
 + |x-x_0| < \beta\, ;\, \bar K(T- t )
 \leq \beta/2 \}.
$$
In fact this property can be checked by an immediate computation
noticing that, for $x = x_0$, $\bar K (T-t) + |x-x_0| = \bar K (T-t)
\leq \beta/2$ and therefore all the derivatives of $\psi$ are equal to
0 at that point. Finally, for any $\kappa$, set 
$$\chi (t,x) = \inf_{x_0}\, \chi_{x_0} (t,x)= \psi( \bar K (T-t) +
 d(x,D))
 + \eta (T-t).$$ 
Then $\chi$ is well-defined and is still a strict supersolution of
 (\ref{timedepeq2}) in the domain 
$${\mathcal C}:=  \{ (t,x) \in \tilde D_{\tau} \ : \ \bar K (T-t)
 + d(x,D) < \beta\, ;\, \bar K(T- t )
 \leq \beta/2 \}.
$$
Now we want to compare $w$ and $\chi$ in ${\mathcal C}$. From the
definition of $\chi $, it is clear that it tends to infinity on the
boundary of the domain ${\mathcal C}  $. Using
 this remark together with standard comparison arguments yield that
 $w \leq \chi$ in ${\mathcal C}$.
For any $x \in D$,
$$w( t, x) \leq \chi( t, x) =  \psi( \bar K (T-t)) +\eta (T- t) = \eta (T- t)
$$
since  $\bar K (T- t) \leq \beta/2$.
Letting $\eta$ tend
 to $0$, we deduce that for $x\in D$ and $\tau \leq t \leq T$ 
$$w(x,t) \leq  \ 0 .$$
Hence $\tilde m (\tau) \leq 0 $ 
 for all $T_0 \leq \tau < T $. Thus 
$\tilde m (T_0) \leq 0 $
and the proof of the first step is complete.\\

\noindent{\bf Proof of Step 2 :} By the definition of $M$ and $N$, and in particular the fact that $g$ takes values in $D$, we deduce that, if $x \in A$, then for any $v_0$
$$ \mu M  u_1(t, x) \leq \mu u_1(t, g (x,v_0)) + \mu C_a (x_\kappa, v_0) \leq u_2(t, g (x,v_0)) + C_a (x_\kappa, v_0) + (\mu-1)C'\; ,$$
and therefore, taking the infimum on $v_0$, $\mu M  u_1(t, x)\leq M
u_2 (t, x) +(\mu-1)C'$. A similar property holds for $Nu_1$ and $N
u_2$. For the boundary points of A and C, we have to argue on the same
lines as in Lemma 4.2, to show that the correct inequalities hold.

We immediately deduce that $w(t,x) \leq 0 $ for $T_0 \leq t \leq T$
and $x \in A$ because 
$\mu  u_1(t, x) - u_2 (t,x) \leq \mu M  u_1(t, x)- M  u_2 (t, x) \leq 0.$\\

\noindent{\bf Proof of Step 3 :} It remains to show that $w \leq 0$ in $[T_0,T]
\times (\Om \setminus A)$. To do so, we first prove that $w$ is a
subsolution of (\ref{timedepeq1}) in the set $\{ w>0 \}$. We only
sketch this proof since it follows arguments we already used above.

The only difficulty comes from the points $(t,x)$ in $C$ since we may have the inequality $u_2 (x,t) \geq N u_2 (x,t)$ (and not the inequality associated to the HJB equation). But, if this is the case, then, by Step~2,
$$ \mu u_1 (t,x) - u_2 (t,x) \leq \mu N  u_1(t, x)- N  u_2 (t, x) \leq (\mu-1)C' < 0\; .$$
Therefore this case cannot happen on the set $\{ w>0 \}$ and the same
arguments
 as in Step 1 above show that $w$ is a subsolution of (\ref{timedepeq1}) in the set $\{ w>0 \}$.

To conclude, we have to modify the comparison argument by building ``friendly giant functions'' as the function $\chi$ above but in a slightly different way in order to take into account the $|x|$-dependence in a suitable way. This also allows us to use a ``local'' comparison argument.

For any interior point $x_0 \in \Om_\kappa  \setminus  A_\kappa$ and
for $r,c >0$, 
we introduce the set
$$ {\cal{O}} {(x_0, r)} = \{(t, x) \in [T_0, T] \times \Om \, : \, |x-x_0| < r \sqrt{c(t - T_0 )} \} \; .$$  
Next we define $\chi : {\cal{O}} {(x_0, r)} \to \R  $ by
$$ \chi (t, x) = \frac {1} { r^2 c( t- T_0 ) - |x-x_0|^2 }\; . $$

We need to show that $\chi$ is a strict supersolution of (\ref{timedepeq1}) in $ {\cal{O}} {(x_0, r)}$. By simple calculations, we get
\beas 
\chi_t &=& \frac {- c  r^2} {( r^2 c( t- T_0 ) - |x - x_0|^2)^2 }
 \\
 D\chi &=& \frac { 2 (x - x_0)} {  (r^2 c( t- T_0 )  - |x-
   x_0|^2 )^{2}   }
 \eeas
Since $|x| \leq |x-x_0| + |x_0|   $ and $|x-x_0|  \leq  r \sqrt{c(T - T_0 )} $, we have
$$
2 F ( 1 + |x|) |x-x_0| 
 \leq 2 F(  1 +  r \sqrt{c(T - T_0 )} + |x_0| ) r\sqrt{c(T - T_0 )} \; .
$$
Finally if we choose $ r \geq r (x_0) = ( 1 + |x_0| ) $ and $c$ such
that 
$\sqrt{c} >  \frac{2 F\sqrt{(T - T_0 )}}{1-2F(T - T_0 )} $, we have
$$2 F ( 1 + |x|) |x-x_0|  < c r^2 .
$$
Note that such a choice of $c$ is possible because of ( \ref{timechoiceeq1}).
Substituting and using above inequality we can conclude that $\chi$ is
a  strict supersolution of (\ref{timedepeq1}) in ${\cal{O}}{(x_0, r)}$.

By definition, $ \chi$ goes to $ \infty$ on the lateral boundary of
${\cal{O}} {(x_0, r)}$ and hence, on that part of the lateral boundary contained inside 
$[T_0, T] \times (\Om_\kappa \setminus  A_\kappa )$, we have no
problem to compare $w$ and $\chi$. Of course, $w \leq \chi$ on
$[T_0, T] \times   A_\kappa$ since  $\chi \geq 0$ and $w\leq 0$ there by step 2. Therefore we can conclude that $w \leq \chi$ in ${\cal{O}} {(x_0, r)}$. 

Now we let $r$ go to infinity to obtain
$$ w \leq  0 \quad \hbox{in  } [T_0, T] \times \Om_\kappa
\setminus  A_\kappa \; ,$$
for every $\kappa  $. The proof is now complete since we can iterate the argument in time according to the fact that $T-T_0$ depends only on $\beta$ and $\bar K$.
\end{proof}

\medskip

Finally we show that the value
function naturally satisfies the condition (\ref{hypid}) under some reasonable
additional assumptions. For that it is enough to show that
$$ | V( t, x)- \tilde  h (x)| \to 0
$$
as $ t$ approaches $T$, uniformly for bounded $ x  \ \in \ \Omega $.

If $t$ is close to $T$ and $|x| \leq r$, it is easy to see that the trajectory $X$ starting from $x$ at time $t$, remains in a bounded subset of $\Omega$ on the time interval $[t,T]$. This is true even after the jumps, since the jump destination set $D$ is bounded. Because of (C2), the number of jumps is finite and estimated uniformly in
$x$ by a fixed number depending on $r$. Therefore, the integral of $K$ gives a $O(T-t)$ whatever the trajectory is.

Next we estimate $V( t, x)- \tilde  h (x)$ from below by considering a
$\varepsilon$-optimal control $(u(\cdot), v, \xi_i,
{X(\xi_i)}^\prime)$ (we drop the 
dependence with respect to $\varepsilon$ for the sake of simplicity of notations); we have
\begin{eqnarray*}
V(t, x) + \varepsilon & \geq &\int\limits_{t}^T K (t, X (t),u (t)) dt 
+ \sum\limits_{s \leq \sigma_i < T } C_a (X (\sigma_i^-) , v )
  \nonumber \\
 &+& \sum\limits_{s \leq \xi_i < T } C_c (X(\xi_i^-), {X(\xi_i)}^\prime) + h( X_x(T))
\end {eqnarray*}
But, using that
$$\tilde h
 (X(s))= M\tilde h (X(s)) \quad \hbox{ if $X(s) \in A$}$$
or 
$$\tilde h (X(s)) = \min  (h(X(s), N h (X(s))) \quad \hbox{ if $X(s) \in C$}, $$
we deduce 
$$ C_a (X (\sigma_i^-) , v ) \geq \tilde h (X (\sigma_i^-)) - \tilde h (X (\sigma_i^+))\; ,$$
$$ C_c (X(\xi_i^-), {X(\xi_i)}^\prime) \geq \tilde h (X(\xi_i^-)) - \tilde h (
{X(\xi_i)}^\prime)\; .$$

Therefore
\begin{eqnarray*}
V(t, x) + \varepsilon & \geq &O(T-t)
+ \sum\limits_{s \leq \sigma_i < T } \left(\tilde h (X (\sigma_i^-)) - 
\tilde h (X (\sigma_i^+))\right)
  \nonumber \\
 &+& \sum\limits_{s \leq \xi_i < T } \left(\tilde h (X(\xi_i^-)) - \tilde h (
{X(\xi_i)}^\prime)\right)+ \tilde h( X_x(T))
\end {eqnarray*}
In order to conclude, it is enough to rearrange the terms of the
right-hand side and to make it appear as a sum of differences of
values $\tilde h$ between two jumps (after one jump and before the
next one). The uniform continuity of $\tilde h$ and the estimates on
the $f$-ode (\ref{1}), lead to the desired property, since $\varepsilon$ is arbitrary.

Conversely, for any $u(\cdot)$, we write
\begin{eqnarray*}
V(t, x) & \leq &\int\limits_{t}^T K (t, X (t),u (t)) dt 
+ \sum\limits_{s \leq \sigma_i < T } C_a (X (\sigma_i^-) , v )
  \nonumber \\
 &+& \sum\limits_{s \leq \xi_i < T } C_c (X(\xi_i^-), {X(\xi_i)}^\prime) + h( X_x(T))
\end {eqnarray*}
but this time, when we touch either $A$ or $C$, we choose the jumps in order to have
$$ C_a (X (\sigma_i^-) , v ) \leq \tilde h (X (\sigma_i^-)) - \tilde h (X (\sigma_i^+)) + \varepsilon \; ,$$
$$ C_c (X(\xi_i^-), {X(\xi_i)}^\prime) \leq \tilde h (X(\xi_i^-)) - \tilde h (
{X(\xi_i)}^\prime) + \varepsilon \; .$$
By the same arguments as above, we again obtain
$$ | V( t, x)- \tilde  h (x)| \leq \hat m(T-t) \quad \hbox{ for any }|x| \leq r $$
where $\hat m(\tau) \to 0$ as $\tau \to 0$ depend only on the $L^\infty$ norm of $K$ and on the modulus of continuity of $\tilde h$ on the bounded subset where the trajectory $X$ lives.

The above inequalities show that $V$ assumes the terminal data $\tilde h$ uniformly w.r.t $i$, which is exactly  (\ref{hypid}).


\begin{thebibliography}{99}

\bibitem{bak} {\sc A.Back, J.Gukenheimer  and M. Myers},  {\em A dynamical simulation facility for hybrid systems},
 Hybrid Systems, Lect. Notes in Comp Sci. vol 736, R.L. Grossman, A. Nerode, A.P.Rava and H. Rischel Eds.  Springer, New York (1993) 

\bibitem{ba1} {\sc G. Barles}, {\em Deterministic Impulse Control
  Problems}, SIAM Jl on control and Optimization, 23, (1985), pp 419-432


\bibitem{ba94} {\sc G. Barles}, {\em Solutions de viscosite  des equations de Hamilton Jacobi,}
 vol 17, of Matematiques et Applications. Springer, Paris, (1994)

\bibitem{ba03} {\sc G. Barles, S. Biton and O. Ley}, {\em Uniqueness for Parabolic 
 equations without growth condition and applications to the mean curvature flow in $\R^2$ },
 Jl of Differential Equations, 187, (2003), pp 456-472

\bibitem{bardi} { \sc M. Bardi  and Capuzzo Dolcetta},  {\em Optimal 
 Control and Viscosity solutions of Hamilton-Jacobi-Bellman
 equations,} Birkhauser, Boston (1997)


 \bibitem{bor} {\sc M. S. Branicky, V. Borkar and S. Mitter}, {\em A 
Unified Framework for Hybrid Control problem,} IEEE Transactions on 
 Automated control, 43 (Jan 1998), pp 31-45 

\bibitem{br}{ \sc M.S. Branicky},   {\em Studies in hybrid systems: 
Modeling, analysis and control,} Ph.D.dissertation, 
Dept. Elec. Eng. and Computer sci., MIT Cambridge (1995)


\bibitem{CIL}
M.G Crandall, H.Ishii and P.L Lions. 
\newblock User's guide to
viscosity solutions of second order partial differential
equations.
\newblock Bull. Amer. Soc. {\bf 27} (1992), pp 1-67.


\bibitem {SDMR1} {\sc S. Dharmatti and M. Ramaswamy } {\em Hybrid Control System and viscosity Solutions,} 
SIAM Jl  on Control and Optimization 34 (2005), pp 1259-1288

\bibitem {SDMR2} {\sc S. Dharmatti and M. Ramaswamy} {\em  Zero Sum Differential Games Involving Hybrid Controls,}  Jl
 of Optimization Theory and Applications 128 (2006), pp 75-102

\bibitem{OL}
Ley, O.
\newblock Lower-bound gradient estimates for first-order {H}amilton-{J}acobi
  equations and applications to the regularity of propagating fronts.
\newblock {\em Adv. Differential Equations}, 6(5):547--576, 2001.


\bibitem{var} {\sc P. P. Varaiya},  {\em Smart cars on smart roads: 
problems of control,} IEEE Transactions on Automated Control, 38 (1993), 
pp 195-207 



\end{thebibliography}
\end{document}